\def\misajour{04/02/2020} 
\documentclass[11pt,leqno]{article} 
 
\usepackage[applemac]{inputenc}
\usepackage[pdftex]{graphicx}
\usepackage{apalike}
 \usepackage%[backref]%,
{hyperref} 

\usepackage{xy}
\usepackage{xypic} 
 
\usepackage{multicol}
\usepackage{xcolor}
\usepackage{amsmath}
\usepackage{amssymb}
\usepackage{amsthm} 
\def\C{{\mathbb C}}
\def\K{{\mathbb K}}
\def\N{{\mathbb N}}
\def\Q{{\mathbb Q}}
\def\R{{\mathbb R}}
\def\Z{{\mathbb Z}}

\def\rmd{{\mathrm d}}
\def\rme{{\mathrm e}}

\def\and{ \hbox{ and } }

\def\phitilde{\tilde\phi}
\def\gammabar{\overline{\gamma}}

\newtheorem{theorem}{Theorem}

\def\adots{\mathinner{\mkern2mu\raise1pt\hbox{.}
\mkern3mu\raise4pt\hbox{.}\mkern1mu\raise7pt\hbox{.}}}

\begin{document}
 
 \hsize 16 true cm 
\vsize 20 true cm 
\hoffset -1 true cm 
 
 \hfill 
 
 \null
 \vskip -3 true cm
 
 \noindent
 31th meeting of the Journées Arithmétiques 
 \hfill
 \emph{\misajour}
 
 \begin{center}
 \LARGE
 \bf 
 
 \medskip
 Integer--valued functions, Hurwitz functions 
 
 and related topics: a survey
 
 \large
 \medskip
 by
 
 \medskip
 \sc
 Michel Waldschmidt
 
 \end{center}

 \section*{Abstract}
 An integer--valued function is an entire function which maps the nonnegative integers $\N$ to the integers. An example is $2^z$.
 A Hurwitz function is an entire function having all derivatives taking integer values at $0$. An example is $\rme^z$. 
 
 Lower bound for the growth order of such functions have a rich history. Many variants have been considered: for instance, assuming that the first $k$ derivatives at the integers are integers, or assuming that the derivatives at $k$ points are integers. These as well as and many other variants have been considered. We survey some of them.

\tableofcontents
 
 \section{On the sequence of values $(f(0), f(1), f(2),\dots)$}
 
 The topic of integer-valued entire functions was initiated by P\'olya's fundamental result on transcendental entire functions taking their values in $\Z$ at each point in $\{0,1,2,\dots\}$. 
 
 \subsection{Order and type of an entire function}
 
Let $f$ be an entire function. For $r\ge 0$ we define $|f|_r=\sup_{|z|\le r|}|f(z)|$. 
From the maximum modulus principle we deduce $|f|_r=\sup_{|z|=r|}|f(z)|$. 

The order of an entire function $f$ is 
$$
\varrho(f)=\limsup_{r\to\infty} \frac{\log\log |f|_r}{\log r},
$$
while the exponential type of an entire function is 
$$
\tau(f)= \limsup_{r\to\infty} \frac{\log |f|_r}{r}=\limsup_{n\to\infty} |f^{(n)}(z_0)|^{1/n} \quad (z_0\in\C).
$$ 
We use the notation: 
$$
f^{(n)}(z)= 
\frac{\rmd^n}{\rmd z^n}
f(z).
$$%
If the exponential type is finite, then $f$ has order $\le 1$. If $f$ has order $<1$, then the exponential type is $0$. Tor $\tau\in\C\setminus\{0\}$, the function $\rme^{\tau z}$ has order $1$ and exponential type $|\tau|$. 

For $\varrho>0$, we define 
$$
\tau_\varrho(f)= \limsup_{r\to\infty} \frac{\log |f|_r}{r^\varrho}, 
$$ 
so that $\tau_1(f)=\tau(f)$. If $f$ has order $<\varrho$, then $\tau_\varrho(f)=0$. If $f$ has order $>\varrho$, then $\tau_\varrho(f)=+\infty$. 

\subsection{Functions vanishing at $0,1,2,\dots$}

Let us consider the entire functions which vanish at each point in $\N=\{0,1,2,\dots\}$. An example of such a function is the Weierstrass canonical product for $\N$, namely \cite[\S~4.41 and \S~8.4.(vi)]{MR3728294}, %E.C. Titchmarsh, The Theory of Functions (2nd Ed) (1939) Oxford University Press 
\cite[Chap.~6]{Shorey}
$$
\frac{1}{\Gamma(-z)}=
 -z\rme^{-\gamma z}\prod_{m=1}^\infty\left(1-\frac z m\right)\rme^{z/m}
 $$
 (Hadamard product for the Gamma function). This is an entire function of order $1$ and infinite type
 \cite[Chap.~2 \S~1.3 (17)]{zbMATH03416363}. % Gelfond calcul des différences finies. 1952
The ``smallest'' entire functions vanishing at each point in $\N$ is 
$$
\sin(\pi z)=\pi z \prod_{m\in\Z\setminus\{0\}} \left(1-\frac{z}{m}\right) \rme^{z/m}
$$ 
which in fact vanishes at each point in $\Z$
(Weierstrass canonical product for $\Z$, Hadamard product for the sine function), which has order $1$ and exponential type $\pi$.
See \cite[\S~3.23 and \S~8.4.(v)]{MR3728294}, %E.C. Titchmarsh, The Theory of Functions (2nd Ed) (1939) Oxford University Press 
and 
\cite[\S~9.4]{MR0068627}. % Boas 1954. 
Indeed, 
 %Fritz David Carlson. 
 a well--known theorem from the thesis of F.~Carlson in 1914 using the Phragmén--Lindelöf principle 
states that there is no nonzero entire function $f$ of exponential type $<\pi$ satisfying $f(\N)=\{0\}$
-- see for instance 
\cite[\S~5.81]{MR3728294}, %E.C. Titchmarsh, The Theory of Functions (2nd Ed) (1939) Oxford University Press 
\cite[Chap.~2, \S~3 Th.~VII and Chap.~3, \S~3.1]{zbMATH03416363}, % Gelfond calcul des différences finies. 1952
\cite[Th.~9.2.1]{MR0068627}, % Boas 1954. 
\cite[Chap.~IV \S 19]{MR0162914}% Boas Buck 1964
.
\subsection{Integer--valued functions}

An \emph{integer--valued function} is an entire function $f$ such that $f(n)\in\Z$ for all $n\ge 0$. A basis of the $\Z$--module of polynomials in $\C[z]$ which map $\N$ to $\Z$ is given by the polynomials 
\begin{equation}\label{Equation:BinomialPolynomials}
1,z,\frac {z (z-1)} 2, \dots,\frac {z (z-1)\cdots (z-m+1)} {m!},\dots
\end{equation}
\cite{MR1367962}, % Narkiewicz 1995
\cite{MR1421321}, % Cahen Chabert 1997
which are sometimes called (in papers from transcendental number theory) the Feld'man's polynomials, after N.I.~Feldman introduced them for improving Baker's lower bounds for linear forms in logarithms of algebraic numbers -- see for instance \cite{zbMATH06841731}. %Livre Bugeaud

The simplest example of a transcendental integer--valued function is $2^z$. 
In 1915, 
in his seminal paper
\cite{Polya1915}, 
 G.~P\' olya proved that if $f$ is an entire function such that
$$
\lim_{r\to\infty}|f|_r 2^{-r}\sqrt r=0,
$$
then $f$ is a polynomial. 
 
The method of P\' olya is to expand $f$ as an interpolation series 
\begin{equation}\label{Equation:InterpolationEntiers}
f(z)=a_0+a_1z+a_2\frac {z (z-1)} 2+\cdots+a_m\frac {z (z-1)\cdots (z-m+1)} {m!}+\cdots
\end{equation}
and to estimate the coefficients $a_m$. See for instance 
\cite[Th.~9.12.1]{MR0068627}. % Boas 1954. .
 
P\'olya conjectured that the conclusion of his theorem should be true under the weaker assumption 
$$
\lim_{r\to\infty}|f|_r 2^{-r} =0.
$$
This refinement was achieved by 
\cite{zbMATH02611157}. %G.H.~Hardy in 1917 
See also \cite{zbMATH02606445} % Landau 
and
\cite[\S~11 Th.~11]{zbMATH02532117}. %Whittaker 1935
In 
\cite{zbMATH02606446}, % Polya 1920 
the same conclusion is obtained under the weaker assumption 
$$
\limsup_{r\to\infty}|f|_r 2^{-r} <1; 
$$
further, if $|f|_r=O(2^r r^k)$ for some $k>0$, then $f(z)$ is of the form $P(z)2^z+Q(z)$ where $P$ and $Q$ are polynomials with rational coefficients. 
It is interesting to note that there are only countably many such functions. 

Twenty years later, \cite{zbMATH02502712} % A.~Selberg in 1941
proved that if an integer--valued function $f$ satisfies 
$$
\tau(f)\le \log 2+\frac{1}{1500},
$$
then $f(z)$ is of the form $P(z)2^z+Q(z)$. 
\cite{zbMATH02500466} % Ch.~Pisot In 1942, 
went one step further. 
If an integer--valued function $f$ has exponential type $\le 0.8$, then $f$ is of the form
$$
P_0(z)+2^z P_1(z)+\gamma^z P_2(z)+\gammabar^z P_3(z),
$$ 
where $P_0,P_1,P_2,P_3$ are polynomials and 
$$
\gamma= 
\frac{3+i\sqrt 3}{2},
\quad
\gammabar = 
\frac{3-i\sqrt 3}{2}
$$
are the roots of the polynomial $z^2-3z+3$. 
See \cite[Th.~9.12.2]{MR0068627}, % Boas 1954. 
\cite{MR29984}. % Buck 1948
This contains the result of Selberg, since $|\log \gamma|=0.758\dots>\log 2=0.693\dots$, 
Pisot 
\cite{MR16483} % CRAS Ch.~Pisot In 1946, 
obtained more general results for functions of exponential type $<0.9934\dots$
with additional terms; he also investigated the growth of transcendental entire functions $f$ with values $f(n)$ close to integers: $f(n)=u_n+\epsilon_n$ with $u_n\in\Z$, $|\epsilon_n|<\kappa^n$ for sufficiently large $n$ and $0<\kappa<1$. Functions which are almost integer-valued were already considered by
\cite{zbMATH02573086} % Gel'fond 1929
and by
\cite[Th.~6.3]{MR19118}. % Buck 1946
Besides, 
\cite[Th.~7.1]{MR19118} % Buck 1946
considered the problem of transcendental entire functions taking prime values at the integers. 

\subsection{Completely integer--valued functions}

A \emph{completely integer--valued function} is an entire function which takes values in $\Z$ at all points in $\Z$. 

Let $u>1$ be a quadratic unit, root of a polynomial $X^2+aX+ 1$ for some $a\in\Z$. 
Then the functions 
$$
u^z+u^{-z}
\and
\frac{u^z-u^{-z}}{u-u^{-1}}
$$%
are completely integer--valued functions of exponential type $\log u$.

Examples of such quadratic units are the roots of the polynomial $X^2-3X+1$:
$$
\theta=\frac{3+\sqrt 5} 2, \quad
\theta^{-1}=\frac{3-\sqrt 5} 2\cdotp
$$%
Hence, examples of completely integer--valued functions are 
$$
\theta^z+
\theta^{-z}
\quad\hbox{and}\quad 
\frac{1}{\sqrt 5} (\theta^z-\theta^{-z}),
 $$
 both of exponential type $\log \theta=0.962\, 423\dots$
 
 Notice that $\theta=\phi^2$ where $\phi$ is the Golden ratio $\frac{1+\sqrt 5} 2$. Let $\phitilde=-\phi^{-1}$, so that 
$$
X^2-X-1=(X-\phi)(X-\phitilde). 
$$% 
For any $m\in\Z$ we have 
$$
\phi^m+\phitilde^m\in\Z.
$$%
The function $\phi^z=\exp(z\log \phi)$ has exponential type $\log \phi<\log 2 $ and we have $
\log|\phitilde|=-\log \phi$. 
However, $\phi^z+\phitilde^z$ not a counterexample to P\'olya's result on the growth of transcendental integer-valued entire functions: indeed, while $\phi^z=\exp(z\log \phi)$ is well defined since $\phi>0$, the definition of $\phitilde^z$ requires to choose a logarithm of the negative number $\phitilde=-\phi^{-1}$. With
$\log \phitilde=-\log \phi+i\pi$, the function $\phitilde^z=\exp(z\log\phitilde)$ has exponential type $((\log \phi)^2+\pi^2)^{1/2}=3.178\, 23\dots>\log 2$. 

According to
\cite{Polya1915}, a completely integer--valued function $f$ which satisfies 
$$
\lim_{r\to\infty} |f|_r\theta^{-r}r^{3/2}=0
$$
 is a polynomial. 
 In 
\cite{zbMATH02606446}, % Polya 1920 
it is proved that a completely integer--valued function $f$ which satisfies 
$$
\limsup_{r\to\infty} |f|_r\theta^{-r}r^k<\infty
$$
for some $k>0$ is of the form 
$$
P_0(z)+P_1(z)\theta^z+ P_2(z)\theta^{-z}
$$
where $P_0,P_1,P_2$ are polynomials. 

\cite{zbMATH02601430} % In 1921, F.~Carlson 
gave refined results 
%proved that if a completely integer--valued function satisfies
%$$
%\tau(f)<\log \theta=0.962\, 423\dots,
%$$
%then $f$ is a polynomial. 
%Carlson used 
by means of Laplace transform -- see 
\cite{MR1545027}, %Polya 1929
\cite[\S~10]{zbMATH02532117}, %Whittaker 1935
\cite[Th.~5.2, Cor.~2]{MR29984}, % Buck 1948 
\cite[Chap.~3, \S~2]{zbMATH03416363}, % Gelfond calcul des différences finies. 1952
\cite[\S~9.12]{MR0068627}, % Boas 1954. 
\cite{MR274762}. % Robinson, Raphael M. 1971 . 
According to \cite{zbMATH02502712}% A.~Selberg in 1941
, if a completely integer--valued function satisfies
$$
\tau(f)\le \log \theta +2\cdot 10^{-6},
$$
then $f$ is of the form
$$
P_0(z)+P_1(z)\theta^z+ P_2(z)\theta^{-z}
$$
where $P_0,P_1,P_2$ are polynomials. 

\subsection{Further results}
 
Among the surveys on these topics, let us quote
\cite{zbMATH03623678}, %by F.~Gramain 1978
\cite{MR920559}, %by F.~Gramain 1985
\cite{MR1050972} %F.~Gramain F.J.~Schnitzer 1989
and
\cite{MR1044109}. % Gramain 1988 

Further results have been proved on the growth of entire functions satisfying $f(\N)\subset \K$ or $f(\Z)\subset \K$, where $\K$ is the field of algebraic numbers. 
Assuming suitable growth conditions on $|f|_r$, on the algebraic numbers $f(m)$ and all of their conjugates, one deduces that $f$ is a polynomial. Interpolation formulae yield such results, and also methods from transcendental number theory. We come back to this topic in \S~\ref{SS:Transcendence}. 

 \section{On the sequence of values $(f(0), f'(0), f''(0),\dots)$}. 
 
 In place of the values of $f$ at the integers, we now consider the derivatives of $f$ at $0$. The sequence polynomials 
 \eqref{Equation:BinomialPolynomials} is replaced by the polynomials $z^n/n!$, while the interpolation series \eqref{Equation:InterpolationEntiers} is replaced by Taylor expansion. 

 \subsection{Hurwitz functions} 
 
A \emph{Hurwitz function} is an entire function $f$ such that $f^{(n)}(0)\in\Z$ for all $n\ge 0$.
Polynomials in $\Q[z]$ of the form
$$
\sum_{n=0}^N a_n \frac{z^n}{n!}
$$%
with $a_0,a_1,\dots,a_N$ in $Z$ are Hurwitz functions, while the exponential function 
$$
\rme^z=1+z+\frac {z^2}{2}+\frac {z^3}{6}+\cdots+\frac {z^n}{n!}+\cdots
$$
is a transcendental Hurwitz function. A basis of the $\Z$--module of polynomials in $\C[z]$ which are Hurwitz functions is given by $1,z,z^2/2,\dots,z^n/n!,\dots$. 

The first lower bound for the growth of a transcendental Hurwiz function is due to 
 \cite{zbMATH02611160}, % S.~Kakeya (1916) 
who proved that a Hurwitz function satisfying 
$$
\limsup _{r\to \infty} |f|_r \rme^{-r} \sqrt r=0
$$
is a polynomial. This was refined by 
\cite{zbMATH02601316} -- % G.~P\' olya in 1921
see also \cite[Part VIII, Chap. 3, \S~6, n$^\circ$187]{MR1492448}: % Polya Szegö
a Hurwitz function satisfying 
$$
\limsup _{r\to \infty} |f|_r \rme^{-r} \sqrt r<\frac{1}{\sqrt{2\pi}}
$$
is a polynomial. This is best possible for uncountably many functions, as shown by the functions
$$
f(z)=\sum_{n\ge 0} \frac{e_n }{2^n!} z^{2^n}
$$
with $e_n\in\{1,-1\}$ which satisfy 
$$
\limsup _{r\to \infty}  |f|_r \rme^{-r}\sqrt r=\frac{1}{\sqrt{2\pi}}\cdotp
$$

 \subsection{Refined estimates by Sato and Straus} 
More precise results are achieved by D.~Sato and E.G.~Straus thanks to a careful study of the function 
 $$
 \phi(r)=\max_{n\ge 0} \frac{r^n}{n!}\cdotp
 $$
In \cite[Corollary 1 p.~304]{MR0159945} % Sato Straus: rate of growth
and 
\cite[Corollary p.~20]{MR0192030} % Sato Straus: on the rate of growth 
(see also \cite{MR0280715}), %Sato on the type
they proved that for every $\epsilon>0$, there exists a transcendental Hurwitz function with 
$$
\limsup_{r\to\infty} |f|_r  \rme^{-r} \sqrt{2\pi r}\left(1+\frac{1+\epsilon}{24 r}\right)^{-1}  <1,
$$
while every Hurwitz function for which 
$$
\limsup_{r\to\infty}  |f|_r\rme^{-r} \sqrt{2\pi r} \left(1+\frac{1-\epsilon}{24 r}\right)^{-1}  \le 1
$$
is a polynomial.

 \section{Several points and / or several derivatives}
 
 \subsection{Introduction}
 There are several natural ways to mix integer--valued functions and Hurwitz functions: one may include finitely may derivatives in the study of integer--valued functions, yielding to the study of $k$--times integer--valued functions. One may consider entire functions having all their derivatives at $k$ points taking integer values, yielding to the study of $k$--point Hurwitz functions. One may consider entire functions $f$ which satisfy $f^{(n)}(m)\in\Z$ for all $n\ge 0$ and $m\in\Z$, which are the so--called utterly integer--valued functions. Finally, the study of entire functions $f$ such that $f^{(n)}(n)\in\Z$ is related with Abel's interpolation. 
  
 Let us display horizontally the rational integers and vertically the derivatives. 
 
 \goodbreak
 
\begin{center}
\bf
Integer--valued functions: 
horizontal 
\end{center}
$$
 \begin{matrix}
 f&
 \bullet & \bullet & \bullet &\cdots &\bullet& \cdots
 \\
&0 & 1& 2 &\cdots &m&\cdots
 \end{matrix}
 $$
 
 \goodbreak
\begin{center}
\bf
Hurwitz functions: 
vertical
\end{center}
 $$
 \begin{matrix}
&\vdots 
\\
 f^{(n)} 
& \bullet 
 \\
&\vdots
\\
 f'& \bullet 
 \\
 f & \bullet 
 \\
&0
 \end{matrix}\hfill
 $$ 
 
 \begin{center}
\bf
 $2$--times integer--valued functions
 \end{center}
 $$
 \begin{matrix}
 f'&
 \bullet & \bullet & \bullet &\cdots &\bullet &\cdots
 \\
 f&
 \bullet & \bullet & \bullet &\cdots &\bullet &\cdots
 \\
&0 & 1& 2 &\cdots &m&\cdots
 \end{matrix}
 $$
 
\begin{center}
\bf
 $2$--point Hurwitz functions
 \end{center}
 $$
 \begin{matrix}
\vdots&\vdots&\vdots
\\
 f^{(n)}&
 \bullet & \bullet 
 \\
\vdots&\vdots&\vdots
\\
 f'&
 \bullet & \bullet 
 \\
 f&
 \bullet & \bullet 
 \\
&0 & 1 
 \end{matrix}\hfill
 $$ 
 
 \begin{center}
\bf
Utterly integer--valued 
\\
entire functions
\end{center}
$$
 \begin{matrix}
\vdots&\vdots&\vdots& &\vdots & 
\\
 f^{(n)}&
 \bullet & \bullet &\cdots&\bullet &\cdots
 \\
\vdots&\vdots&\vdots&\ddots&\vdots & 
\\
 f'&
 \bullet & \bullet &\cdots&\bullet &\cdots
 \\
 f&
 \bullet & \bullet &\cdots&\bullet &\cdots
 \\
&0 & 1 &\cdots&m &\cdots
 \end{matrix}\hfill
 $$
 
  \goodbreak
 
\begin{center}
\bf
Abel interpolation
\end{center}
 $$
 \begin{matrix}
 \vdots
 & & & & & \adots
\\
 f^{(n)}&
 & & &\bullet & 
 \\
 \vdots
 & & &\adots& & 
\\
 f'&
 & \bullet & & & 
 \\
 f&
 \bullet & & & & 
 \\
&0 & 1 &\cdots &n & \cdots
 \end{matrix}\hfill
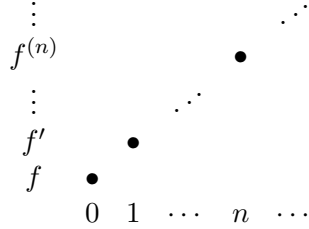
 $$ 
 
 \subsection{$k$--times integer--valued functions}
 The first natural way to mix integer--valued functions and Hurwitz functions is \emph{horizontal}, including finitely may derivatives in the study of integer--valued functions, like in 
 \cite{zbMATH02569599}, % Gel'fond (1929), 
\cite{MR0006569}. %Selberg (1941). 
Let us call \emph{$k$--times integer--valued function} an entire function $f$ such that $f^{(n)}(m)\in\Z$ for all $m\ge 0$ and $n=0,1,\dots,k-1$.
According to 
\cite{zbMATH02573085}, %{\bf Gel'fond} (1929) 
a $k$--times integer--valued function of exponential type $< k\log\left(1 + \rme^{-\tfrac{k-1}{k}}\right)$ is a polynomial. A proof is given in 
\cite{MR0223571}. % Fridman
 Improvements are due to 
 \cite{zbMATH02502712}, % A.~Selberg in 1941
 and
\cite{MR2100701} % Bundschuh Zudilin 2004
who also investigated what could be the best possible results. The best known results so far are due to
\cite{MR2110509}. % Welter 2005

 \subsection{$k$--point Hurwitz functions}\label{SS:kpointsH}
 The second solution is
\emph{vertical}, considering \emph{$k$--point Hurwitz functions}, namely entire functions having all their derivatives at $0,1,\dots,k-1$ taking integer--values. 
This question was first considered by 
\cite{zbMATH03018993}. % Gel'fond 1934
It has been proved by \cite{MR0035822} % Straus (1950)
that the order of such a function is $\ge k$, and this is best possible, as shown by the function 
$\rme^{z(z-1)\cdots(z-k+1)}$. 
 
Precise results are known for 
$k=2$. 
\cite{MR0280715} % Sato on the type
 proved that 
there exist transcendental two point Hurwitz entire functions with 
$$
|f|_r\le \exp\left( r^2+r-\log r+O(1)\right),
$$
while every two point Hurwitz entire functions with 
$$
|f|_r\le C \exp\left( r^2-r-\log r\right)
$$
for some positive constant $C$ must be a polynomial. 

Another example is 
\cite[Th.~3]{MR0035822}: % E. G.~Straus 1950
if $f$ is a transcendental entire function and $f^{(n)}(z)$ is in $\Z$ for $z=0$ and $z=p/q$ with $\gcd(p,q)=1$, then $f$ is at least of order $2$ and type $\tau_2(f)\ge q/p$. 
A similar estimate from
\cite[Th.~2]{MR0035822} % E. G.~Straus 1950
holds when $f^{(n)}(z)$ is in $\Z$ for $z=0,p_1/q_1,\dots,p_{k-1}/q_{k-1}$. 
 
For $k\ge 3$ our knowledge is more limited. 
 \cite{MR0055443} % L. Bieberbach 
stated that if a transcendental entire function $f$ of order $\varrho$ 
is a $k$--point Hurwitz entire function, then either $\varrho>k$, or $\varrho=k$ and the type $\tau_k(f)$ of $f$ satisfies $\tau_k(f)\ge 1$. 
However, as noted by 
\cite{MR0223571} % Fridman 1968
and
 \cite[p.~2]{MR0280715}, %Sato 1971 on the typerom 
 this result is not true. 
 Indeed, the polynomial
$$
a(z)=
\frac 1 2 z(z-1)(z-2)(z-3) 
$$% 
can be written
$$
a(z)= \frac 1 2 z^4-3z^3- \frac {11} 2z^2-3z, 
$$%
hence it satisfies $a'(z)\in\Z[z]$;
it follows that the function $\rme^{a(z)}$
is a transcendental $4$-point Hurwitz function of order $\varrho=4$ and $\tau_4(f)=1/2$. 

\subsection{Utterly integer--valued functions}
An \emph{utterly integer--valued function} is an entire function $f$ which satisfies $f^{(n)}(m)\in\Z$ for all $n\ge 0$ and $m\in\Z$. 
From either the results on $k$--point Hurwitz functions or the results on $k$--times integer--valued functions, it follows that a transcendental utterly integer--valued function must be of infinite order. One also deduces the irrationality of any nonzero period of a nonconstant Hurwitz function of finite order
-- compare with \S~\ref{SS:Transcendence}. 

 \cite{MR0040481} % Straus (1951): 
suggested that transcendental utterly integer--valued functions may not exist.
\cite{MR0223571} % Fridman 1968
showed that there exists transcendental utterly integer--valued function $f$ with 
 $$
\limsup_{r\to\infty}\frac{\log\log |f|_r}{r} \le \pi
$$
and proved that a transcendental utterly integer--valued function $f$ satisfying
$$
\limsup_{r\to\infty}\frac{\log\log |f|_r}{r} \ge \log(1+1/\rme)
$$
is a polynomial. 
The bound $ \log(1+1/\rme)$ was improved by 
\cite{MR2146856} % M.~Welter
to $\log 2$. Hence a transcendental utterly integer--valued function grows at least like the double exponential $\rme^{2^z}$. 

\cite{MR789191} % D.~Sato Utterly integer-valued entire functions
constructed a nondenumerable set of utterly integer--valued functions. 
He selected inductively the coefficients $a_n$ with 
$$
\frac{1}{n!(2\pi)^n}\le |a_n|\le \frac{3}{n!(2\pi)^n}
$$
and defined
$$
f(z)=\sum_{n\ge 0} a_n \sin^n(2\pi z).
$$

\subsection{Abel interpolation}

 There is also a \emph{diagonal} way of mixing the questions of integer--valued functions and Hurwitz functions by considering entire functions $f$ such that $f^{(n)}(n)\in\Z$.

The source of this question goes back to \cite{zbMATH06075652}, % Abel: Oeuvres
 \cite{abel1881}. % référence précise 1881
 The related interpolation problem was studied by
\cite{zbMATH02706584}, % Halphén CRAS 1882
\cite{zbMATH02705364}. % Halphén BSMF 1882
See also 
\cite{zbMATH02685243}, % Vincente Pareto 1892
\cite{zbMATH02563461}, %Gontcharoff 1930
\cite{MR19118} % Buck 1946
and
\cite[\S~7 and \S~10 Corollary 5]{MR0029985}. % Buck interpolation series 1948
%and
Further references are 
\cite[Chap.~III]{zbMATH02532117}, %Whittaker 1935
\cite[Chap.~III, \S~3.2]{zbMATH03416363}, % Gelfond calcul des différences finies. 1952
\cite[\S~9.10]{MR0068627}, % Boas 1954,
\cite{MR1154472}. % Bézivin 1992

The sequence of polynomials $(P_n)_{n\ge 0}$ defined by 
 $$
 P_0=1, \quad P_n(z)=\frac 1 {n!} z(z-n)^{n-1} \qquad (n\ge 1) 
 $$
 was introduced by 
 \cite{abel1881}. % référence précise
 These polynomials satisfy 
 $$
 P'_n(z)=P_{n-1}(z-1) \quad (n\ge 1). 
 $$
 Therefore $P_n^{(k)}(k)=\delta_{kn}$ for $k$ and $n\ge 0$. One deduces that any polynomial $f$ has a finite expansion
 $$
 f(z)=\sum_{n\ge 0} f^{(n)}(n)P_n(z).
 $$
It was proved by 
\cite{zbMATH02705364} % Halphén BSMF
(see also
\cite[p.~31]{zbMATH02563461}) %Gontcharoff
that such an expansion (with a series in the right hand side which is absolutely and uniformly convergent on any compact of $\C$) holds also for any entire function $f$ of finite exponential type $<\omega$, where $\omega=0.278\, 464\, 542\dots$ is the positive real number defined by $\omega \rme^{\omega+1}=1$.
 The proof rests on Laplace transform -- see 
\cite[Chap.~3 \S~2 p.209]{zbMATH03416363}. % Gelfond calcul des différences finies. 1952
 
Lower bounds for the growth of entire functions satisfying $f^{(n)}(n)\in\Z$ were investigated by 
\cite{MR95939}. % Bertrandias 1958
See also
\cite{MR244483}. %Wallisser 1969
The method arises from 
\cite{MR16483} %Pisot CRAS
and
\cite{MR16482}. %Pisot CRAS

For $t\in\C$, the function $f_t(z):=\rme^{tz}$ satisfies the functional differential equation 
$$
f'(z)=t\rme^tf(z-1)
$$
with the initial condition $f_t(0)=1$, hence $f_t^{(n)}(n)=(t\rme^t)^n$ for all $n\ge 0$. 
 Let $\tau_0=0.567\, 143\, 290\dots$ 
be the positive real number defined by $\tau_0 \rme^{\tau_0}=1$. We have $f_{\tau_0}^{(n)}(n)=1$ for all $n\ge 0$. 

 Following \cite{MR95939}, % Bertrandias 1958
an entire function $f$ of exponential type $<\tau_0$ such that $f^{(n)}(n)\in\Z$ for all sufficiently large integers $n\ge 0$ is a polynomial. 
More precisely, let 
$\tau_ 1$ 
be the complex number defined by $\tau_1 \rme^{\tau_1}=(1+i\sqrt 3)/2$; its modulus is $|\tau_1|=0.616\dots$. Then an entire function $f$ of exponential type $<|\tau_1| $ such that $f^{(n)}(n)\in\Z$ for all sufficiently large integers $n\ge 0$ is of the form $P(z)+Q(z) \rme^{\tau_0z}$, where $P$ and $Q$ are polynomials. 
 
\section{Variants}

We do claim to quote all variations and related works around this theme; we only give a few examples among those which would deserve further surveys. 
\\
$\bullet$ 
$q$ analogues and multiplicative versions (geometric progressions): 
\cite{zbMATH02546306}, %Gel'fond 1933
\cite[Chap.~2, \S~3.4, Th.~VIII]{zbMATH03416363}, % Gelfond calcul des différences finies. 1952
\cite{MR215795}, % Gel'fond 1967
\cite{MR0333196}, % Kaz'min A certain problem of A.O. Gel'fond 1973
\cite{MR775551}, %Bézivin 1984
\cite{zbMATH03959890}, % Wallisser 1985
\cite{MR1044109}, % Gramain 1990
\cite{MR1198484}, % Bundschuh 1991
\cite{zbMATH00001020}, % Bézivin 1990
\cite{MR1302505}, % Bézivin 1994 
\cite{MR1313352}, % Bézivin 1994 
\cite{MR1336220}, % Bundschuh Shiokawa 1995
 \cite{MR1794332}, % Welter 2000
\cite{MR2146856}, % Welter 2005
\cite{MR2110509}, % Welter 2005
\cite{MR3204737} %Bézivin 2014 
\dots 
\\
$\bullet$ extensions of the additive and multiplicative versions considered in 
\cite{MR1683633} % Pila, Jonathan and Rodriguez Villegas, Fernando 1999
who introduce the notion of \emph{concordant sequences}; further developments along these lines are in 
\cite{MR618350}, %Perelli Zannier %Perelli, A.; Zannier, U. 
\cite{zbMATH01239698}, % Bézivin 1998
\cite{zbMATH02216976}, % Pila 2002 
\cite{MR1972183}, %Pila 2003
\cite{MR2210649}, % Pila 2005
\cite{MR2429036}, % Pila 
and
\cite{zbMATH05556407}. % Pila 2009
\\
$\bullet$ analogs in finite characteristic:
\cite{MR2837135}, % David Adam
\cite{MR3323345} % Adam Welter
\dots 
\\
$\bullet$ several variables: 
\cite{MR0190092}, % S.~Lang(1965)
\cite{MR185140}, %Gross (1965)
\cite{MR0213304}, %A.~Baker (1967)
\cite{MR0396988}, %V.~Avanissian and R.~Gay \ (1975)
\cite{MR0427671}, % F.~Gramain 1977
\cite{zbMATH03623678}, % F.~Gramain 1978
\cite{MR552466}, % Bundschuh 1980
\cite{MR775551}, %Bézivin 1984
 \cite{MR920559}, % F.~Gramain 1986
 \cite{MR1044109} % F.~Gramain 1990, 
\dots 
 \\
 $\bullet$ definable functions and minimal models: 
\cite{MR3007661}.% Jones Thomas Wilkie 2012

%% p-adic functions
% V. Laohakosol, J. H. Loxton et A. J. van der Poorten, Integer-valued p-adic functions, in Number Theory, Vol. II, (Budapest/Hung. 1987), Colloq. Math. Soc. János Bolyai 51, 1990, pp. 829–849. MR1058247
% voir le résumé critique dans MR
%++++++++++++++
\section{Connection with transcendental number theory}\label{SS:TNT}
 
 \subsection{From $\Z$ to $\Z[i]$}
 
Historically, the very first steps towards a solution of Hilbert's seventh on the transcendence of $a^b$ came from the extension by 
\cite{zbMATH02576736} %S.~Fukasawa (1928)
and
 \cite{zbMATH02569599} %A.O.~Gel'fond(1929)
of P\'olya's result on integer--valued functions on $\N$ to integer--valued functions on $\Z[i]$. 

According to A.O.~Gel'fond, an entire function $f$ which is not a polynomial and satisfies $f(a+ib)\in\Z[i]$ for all $a+ib\in\Z[i]$ satisfies 
 $$
\limsup_{R\rightarrow\infty}\frac{1}{R^2}\log|f|_R\ge\gamma.
$$
For the proof, Gel'fond expands $f(z)$ into a Newton interpolation series at the Gaussian integers. 
He obtained $\gamma\ge 10^{-45}$. 

Since the canonical product associated with the lattice $\Z[i]$, namely the Weierstrass sigma function 
$$
\sigma(z)=z \prod_{\omega\in\Z[i]\setminus\{0\}}\left(1-\frac{z}{\omega}\right)
\exp\left(
\frac{z}{\omega}+\frac{z^2}{2\omega^2}
\right),
$$
 is an entire function vanishing on $\Z[i]$ of order $2$ with $\tau_2(\sigma)=\pi/2$
 \cite[Part 2, Chap.~I, \S~13]{MR0011320}, % Hurwitz Courant - cité par Polya-Szegö; voir Pogany pour des résultats plus précis 
 \cite[Part IV, Chap. 1, \S~3, n$^\circ$49]{MR1492448}, % Polya Szegö
 the best (largest) admissible value for $\gamma$ satisfies 
 $$
10^{-45}\le \gamma\le \frac \pi 2\cdotp
$$
The study of entire functions $f$ satisfying 
$f(\Z[i])\subset \Z[i]$ was pursued by 
\cite{MR0244486} % A.F.~Cayford (1969 ) 
and
\cite{MR614318}. %L.~Gruman (1979),
The exact value is $\displaystyle \gamma=\frac \pi {2e}$: 
\cite{MR590974} %D.W.~Masser (1980) 
proved the upper bound 
and
 \cite{MR620681} % F.~Gramain (1981) :
 the lower bound. 
 
A side effect of these works is the introduction of the so--called Masser--Gramain--Weber constant
 \cite{MR590974}, %D.W.~Masser (1980) 
 \cite{MR771043}, % Gramain--Weber
an analog of Euler's constant for $\Z[i]$, % introduced by D.W.~Masser. 
which arises in a $2$--dimensional analogue of Stirling's formula:
$$
\delta=\lim_{n\to\infty} \left( \sum_{k=2}^n (\pi r_k^2)^{-1}-\log n\right),
$$
where $r_k$ is the radius of the smallest disc in $\R^2$ that contains at least $k$ integer lattice points inside it or on its boundary.
 In \cite{MR3008857} %Melquiond, Guillaume and Nowak, W. Georg and Zimmermann, Paul
the first four digits are computed: 
$$
1.819776<\delta<1.819833,
$$
disproving a conjecture of \cite{MR693311}. %Gramain

\subsection{Transcendence}\label{SS:Transcendence}
 
 One of the main motivations of studying this kind of problems arises from transcendental number theory. See for instance the application to the Hermite--Lindemann Theorem in
\cite[Chap.~2, \S~3.4, Th.~IX]{zbMATH03416363}. % Gelfond calcul des différences finies. 1952

%In 1950, E. G.~Straus 
\cite{MR0035822} developed the subject of integer--valued functions in connection with transcendental number theory; he deduces the Hermite--Lindemann Theorem from his \cite[Th.~4]{MR0035822}. %already quoted in \S~\ref{SS:kpointsH}. 
However, as he pointed out in a footnote, since his paper was written, there has appeared the paper
\cite{MR0031498} % Th.~Schneider 
\emph{ which contains an approach very similar}, which \emph{is more profound}, and the \emph{results are much more complete}. Schneider's work gave rise to the so--called Schneider--Lang Criterion %S.~Lang 
\cite{MR0156852}. % (1962)

One may notice that a main assumption in the final version of the Schneider--Lang criterion is that the functions satisfy a differential equation; hence this criterion does not contain 
\cite[Satz III]{MR0031498} % Th.~Schneider 
nor 
\cite[Th.~4]{MR0035822} %In 1950, E. G.~Straus 
where no such condition occurs. 
For instance, one easily deduces from either of these two results \emph{the transcendence of any nonzero period of a nonconstant Hurwitz function of finite order}. The main example of course is the transcendence of $\pi$. 

When there is no assumption on the derivatives, one may apply the transcendence method introduced by Th.~Schneider for the solution of Hilbert's seventh problem -- see for instance
 \cite{MR485720}. % MW Polya's Thm Schneider's method 1978 
Examples are given in 
\cite{MR702192}, % Gramain Mignotte 1983 
 \cite{MR867491}, % Gramain Mignotte W 1986
\cite{zbMATH05210315}, %Rochev 2007
\cite{MR2866521} %Rochev 2011
and
\cite{MR2828569}. % Ably 2011 
A connection with the six exponentials Theorem is introduced in 
\cite{MR2429036}. % Pila 

\section{Lidstone interpolation and generalizations} 

We conclude by stating some new results related with integer--valued functions. The proofs are given in 
\cite{MW-Lidstone2points} and \cite{MW-LidstoneSeveralpoints}. 
 
\subsection{Arithmetic result for Poritsky and Lidstone interpolation} 
 %%%%%%%%%%%%%%%%%%%%%%%%%%%%%%%%%%%%%%%%%
 
\begin{theorem}
Let $s_0,s_1,\dots,s_{m-1}$ be distinct complex numbers and $f$ an entire function of sufficiently small exponential type.
If 
 $$
 f^{(mn)}(s_j) \in \Z
 $$%
for all sufficiently large $n$ and for $0\le j\le m-1$, then $f$ is a polynomial. 
 \end{theorem}
 
 For $m=2$ (Lidstone interpolation), 
 with $f^{(2n)}(s_0)\in\Z$ and $f^{(2n)}(s_1)\in\Z$, the assumption on the exponential type $\tau(f)$ of $f$ is 
 $$
 \tau(f)<\min\{1, \pi/|s_1-s_0|\},
 $$%
and this is best possible, as shown by the functions 
 $$
 f(z)=\frac{\sinh(z- s_1)}{\sinh(s_0- s_1)} 
\quad \hbox{and}\quad
f(z)=\sin\left( \pi\frac{z-s_0}{ s_1-s_0}\right),
$$%
which have exponential type $1$ and $ \pi/|s_1-s_0|$ respectively. 

When $|s_1-s_0|\le 2$, there is a nondenumerable set of entire functions $f$ of exponential type $\le 1$ satisfying 
$f^{(2n)}(s_0)=0$ and $f^{(2n)}(s_1)\in\{-1,0,1\}$ for all $n\ge 0$. 

\subsection{Arithmetic result for Gontcharoff and Whittaker interpolation 
} 
 %%%%%%%%%%%%%%%%%%%%%%%%%%%%%%%%%%%%%%%%%
 
\begin{theorem}
Let $s_0,s_1,\dots,s_{m-1}$ be distinct complex numbers and $f$ an entire function of sufficiently small exponential type.
Assume that for each sufficiently large $n$, one at least of the numbers 
 $$
 f^{(n)}(s_j) \quad j=0,1,\dots,m-1
 $$% 
is in $\Z$. Then $f$ is a polynomial. 
 \end{theorem}
  
 In the case $m=2$ with 
 $f^{(2n+1)}(s_0)\in\Z$ and $f^{(2n)}(s_1)\in\Z$ (Whittaker interpolation), the assumption on the exponential type of $f$ is 
 $$
 \tau(f)<\min\left\{1, \frac{\pi}{2|s_1-s_0|}\right\},
 $$%
and this is best possible, as shown by the functions 
 $$
f(z)= \frac{\cosh(z- s_1)}{\cosh(s_0- s_1)} 
\quad \hbox{and}\quad
f(z)=\cos\left(\frac{\pi}{2}\cdot \frac{z-s_0}{ s_1-s_0}
\right),
$$%
which have exponential type $1$ and $\frac{\pi}{2| s_1-s_0|}$ respectively. 

When $|s_1-s_0|<\log (2+\sqrt 3)=1.316\dots$, there is a nondenumerable set of entire functions $f$ of exponential type $\le 1$ satisfying 
$f^{(2n+1)}(s_0)=0$ and $f^{(2n)}(s_1)\in\{-1,0,1\}$ for all $n\ge 0$. 

\bigskip 
 %%%%%%%%%%%%%%%%%%%%%%%%%%%%%%%%%%%%%%%%%
 As before, let us display horizontally the points (they are no longer assumed to be the consecutive integers) and vertically the derivatives. 
 \hfill\break\null\qquad 
$\bullet$ interpolation values
\qquad
$\circ$ no condition 
 
 \goodbreak
 \begin{center}
\bf
Lidstone interpolation 
\end{center}
 $$
 \begin{matrix}
 \vdots
 & \vdots& \vdots
\\
 f^{(2n+1)}&
 \circ
 & 
 \circ\\
 f^{(2n)}&
 \bullet & \bullet 
 \\
 \vdots
 & \vdots& \vdots
 \\
 f''&
 \bullet & \bullet 
\\
 f'&
 \circ
 & 
 \circ
 \\
 f&
 \bullet & \bullet 
 \\
&s_0 & s_1 
 \end{matrix}\hfill
 $$
 
 \goodbreak
 \begin{center}
\bf
Whittaker interpolation
\end{center}
 $$
 \begin{matrix}
 \vdots
 & \vdots& \vdots
\\
 f^{(2n+1)}& \bullet 
 & \circ
\\
 f^{(2n)}& \circ
& \bullet 
 \\
 \vdots
 & \vdots& \vdots
 \\
 f''& \circ
 & \bullet 
\\
 f'& \bullet
 & \circ
 \\
 f& \circ
 & \bullet 
 \\
&s_0 & s_1 
 \end{matrix}\hfill
 $$
 
% \bigskip
 
 \goodbreak
 \begin{center}
\bf
Poritsky interpolation ($3$ points)
\end{center}

 $$
 \begin{matrix}
 \vdots
 & \vdots& \vdots& \vdots
\\
 f^{(3n+2)}&
 \circ
 & 
 \circ & 
 \circ\\ f^{(3n+1)}&
 \circ
 & 
 \circ & 
 \circ\\
 f^{(3n)}&
 \bullet & \bullet & \bullet 
 \\
 \vdots
 & \vdots& \vdots& \vdots
 \\
 f^{(\mathrm {iv})}
 &
 \circ
 & 
 \circ& 
 \circ
 \\
 f'''&
 \bullet & \bullet & \bullet 
 \\
 f''&
 \circ
 & 
 \circ& 
 \circ
\\
 f'&
 \circ
 & 
 \circ& 
 \circ
 \\
 f&
 \bullet & \bullet & \bullet 
 \\
&s_0 & s_1& s_2 
 \end{matrix}\hfill
 $$

 \goodbreak 
 \begin{center}
\bf
Gontcharoff interpolation ($3$ points) 

An example with a period 
\end{center}

 $$
 \begin{matrix}
 \vdots
 & \vdots& \vdots& \vdots
\\
 f^{(3n+2)}&
 \circ
 & 
 \circ & 
\bullet
\\ 
f^{(3n+1)}&
 \circ
 & 
\bullet & 
 \circ\\
 f^{(3n)}&
 \bullet & \circ & \circ
 \\
 \vdots
 & \vdots& \vdots& \vdots
 \\
 f^{(\mathrm {iv})}
 &
 \circ
 & 
\bullet& 
 \circ
 \\
 f'''&
 \bullet & \circ & \circ
 \\
 f''&
 \circ
 & 
 \circ& 
\bullet
\\
 f'&
 \circ
 & 
\bullet& 
 \circ
 \\
 f&
 \bullet & \circ & \circ 
 \\
&s_0 & s_1& s_2 
 \end{matrix}\hfill
 $$

%\nocite{*}
%\bibliographystyle{plain}
%\bibliographystyle{smfplain}
%\bibliographystyle{apalike}
%\bibliography{Lidstone}

 \vskip 2truecm plus .5truecm minus .5truecm 
 
 \null
 \hfill
\vbox{\hbox{{\sc Michel WALDSCHMIDT,} Sorbonne Université} 
	\hbox{CNRS, Institut Mathématique de Jussieu Paris Rive Gauche} 
	\hbox{75005 Paris, France}
	\hbox{E-mail: \url{michel.waldschmidt@imj-prg.fr}} 
	\hbox{Url: \url{http://www.imj-prg.fr/~michel.waldschmidt}}	
}	

\end{document}